\documentclass{article}

\begin{document}

\title{K\"ahler-Einstein metrics and algebraic geometry}
\author{Simon Donaldson}
\date{\today}
\maketitle

\begin{abstract}
This paper is a survey of some recent developments in the area described by the title, and follows the lines of the author's lecture in the 2015 Harvard Current Developments in Mathematics meeting.   The main focus of the paper is on the Yau conjecture relating the existence of K\"ahler-Einstein metrics on Fano manifolds to K-stability. We discuss four different proofs of this, by different authors,  which have appeared over the past few years. These  involve an interesting  variety of approaches and draw on techniques from different fields.
   \end{abstract}

\newtheorem{thm}{Theorem}
\newtheorem{defn}{Definition}
\newcommand{\dbd}{\partial\overline{\partial}}
\newcommand{\bC}{{\bf C}}
\newcommand{\bP}{{\bf P}}
\newcommand{\db}{\overline{\partial}}
\newcommand{\bR}{{\bf R}}
\newcommand{\Voln}{\beta_{n}}
\newcommand{\Ric}{{\rm Ricci}}
\newcommand{\OmegaEuc}{\Omega_{{\rm Euc}}}
\newcommand{\reg}{{\rm reg}}
\newcommand{\canon}{{\cal KS}(n,\kappa)}
\section{Introduction}

General existence questions involving the Ricci curvature of compact K\"ahler manifolds go back at least to work of Calabi in the 1950's \cite{kn:Cal1}, \cite{kn:Cal2}.  We begin by recalling some very basic notions in K\"ahler geometry.
 \begin{itemize}\item   All the K\"ahler metrics in a given cohomology class can be described in terms of some fixed reference metric $\omega_{0}$ and a potential function, that is  \begin{equation}\omega= \omega_{0}+ i \dbd \phi. \end{equation}
\item A hermitian holomorphic line bundle over a complex manifold has a unique {\it Chern connection} compatible with both structures. A Hermitian metric on the anticanonical line bundle $K_{X}^{-1}=\Lambda^{n}TX$ is the same as a volume form on the manifold. When this volume form is derived from a K\"ahler metric the curvature of the Chern connection can be identified with the Ricci tensor. (In general in this article we will not distinguish between metrics and Ricci tensors regarded as symmetric tensors or $(1,1)$-forms.) \item If the K\"ahler class is a $2\pi$ times an integral class a metric can be regarded as the curvature of the Chern connection on a holomorphic line bundle. The K\"ahler potential parametrising metrics in (1) has a geometrical meaning as a change in the Hermitian metric on the line bundle: $\vert \ \vert =e^{\phi}\vert\ \vert_{0}$.
 \end{itemize}
One of the questions initiated by Calabi was that of prescribing the volume form of a K\"ahler metric in a fixed cohomology class. By the $\dbd$-lemma this is the same as prescribing the Ricci tensor (as a closed $(1,1)$-form in the class $c_{1}(X)$).  This Calabi conjecture was established by Yau in 1976 \cite{kn:Y}. In particular  when $c_{1}(X)=0$ this gives the existence of {\it Calabi-Yau metrics},  with vanishing Ricci curvature. Another question raised by Calabi involved {\it K\"ahler-Einstein metrics} with  \lq\lq cosmological constant'' $\lambda$,  so
$  {\rm Ricci}=\lambda \omega$.  The case $\lambda=0$ is the Calabi-Yau case as above and  and if $\lambda$ is non-zero  we may assume that it is $\pm 1$. The K\"ahler-Einstein condition can then be expressed as saying that we require that the Hermitian metric on the holomorphic line bundle $K_{X}^{-1}$ given by the volume form of $\omega$ realises $\pm \omega$ as the curvature form of its Chern connection.  Explicitly, in complex dimension $n$ the equation to be solved for a K\"ahler potential $\phi$ is
\begin{equation}  (\omega_{0}+ i \dbd\phi)^{n}= \omega_{0}^{n} \exp(\pm \phi +h_{0})\end{equation}
where $h_{0}$ is the solution of $\dbd h_{0}+\omega_{0}=\pm {\rm Ricci}(\omega_{0}) $ given by the $\dbd$-lemma.

In the case when $\lambda=-1$ there is a straightforward existence theorem, established by Aubin and Yau. In fact the proof is significantly simpler than that when $\lambda=0$. The point is that the nonlinearity from the exponential in (2) occurs with a favourable sign.  The condition on the K\"ahler class means that the manifold $X$ is of general type and the metrics are generalisations of metrics of constant curvature $-1$ on Riemann surfaces of genus 2 or more. The harder case is when $\lambda=1$. Then the condition on the K\"ahler class means that $X$ is a {\it Fano manifold}. A result of Matsushima \cite{kn:Mat}, also from the 1950's, shows that existence can fail in this case. Matsushima showed that if there is a solution the holomorphic automorphism group of $X$ is reductive (the complexification of a compact Lie group). This means that Fano manifolds, such as the projective plane blown up in one or two points, with non-reductive automorphism groups cannot support such a metric. The same Calabi-Aubin-Yau scheme which proved existence in the cases $\lambda\leq 0$--the \lq\lq continity method'' (see 4.2 below)--- can be set up in the positive case, but must break down for these manifolds.  More precisely, the difference arises because of the absence of a $C^{0}$-estimate for the K\"ahler potential due to the unfavourable sign in the equation (2).  On the other hand there are many cases where the existence of a solution has been established, using arguments exploiting detailed geometric features of the manifolds and the theory of \lq\lq log canonical thresholds''. The question which arises is to characterise in terms of the complex geometry of the manifold $X$ exactly when a solution exists.  In the early 1990's Yau conjectured that the appropriate criterion should be in terms of the {\it stability} of the manifold $X$ and, after two decades of work by many mathematicians, this is now known to the case. The precise formulation is in terms of a algebro-geometric notion of {\it K-stability}  and the statement is that a Fano manifold admits a K\"ahler-Einstein metric if and only if it is K-stable.

 One of the attractive features of this problem is the range of techniques which can be brought to bear. By its nature, the statement involves an interaction between algebraic geometry and complex differential geometry and, as we shall see,  there are important connections with global Riemannian geometry, with pluripotential theory in complex analysis and with nonlinear PDE. Four different proofs of the main result have appeared up to the time of writing.
\begin{enumerate} \item Deformation of cone singularities. (Chen, Donaldson and Sun\cite{kn:CDS});
\item The continuity method (Datar, Szekelyhidi \cite{kn:DaS});
\item Proof via Kahler-Ricci flow (Chen, Sun, Wang \cite{kn:CSW});
\item Proof by the variational method (Berman, Boucksom, Jonnson \cite{kn:BBJ}-the statement proved here is slightly different).
\end{enumerate}

The purpose of this article is to survey these recent developments.

\section{K-stability}

The notion of \lq\lq stability'', in this context, arose from the study of moduli problems in algebraic geometry and geometric invariant theory.  It is not usually possible to give a good structure to a set of all isomorphism classes of algebro-geometric objects. For example, generic 4-tuples of points in the projective line up to the action of projective transformations are classified by the cross-ratio, but there is no satisfactory way to define the cross-ratio if three or more points coincide. The  idea is that the isomorphism classes of a suitable restricted class of stable objects do form a good space. There is a general circle of ideas relating the algebraic approach to these questions to metric structures and differential geometry involving the notion of a moment map and the \lq\lq equality of symplectic and complex quotients''. In particular there is a large literature, going back to the early 1980's, developing these notions in the framework of gauge theories and results such as the existence of Hermitian Yang-Mills connections on stable vector bundles \cite{kn:UY}. But the author recently wrote a survey which emphasised this side of the story \cite{kn:D2}, so we will not go into it here beyond noting that the K\"ahler-Einstein problem fits  naturally into  wider questions of the existence of constant scalar curvature and extremal K\"ahler metrics; questions which remain largely open.

Typically, stability is defined by a {\it numerical criterion} on {\it degenerations} of the objects in question. In our situation, we consider polarised varieties $(X,L)$, so $L$ is an ample line bundle over $X$ and the sections of $L$ embed $X$ as a projective variety in some projective space. The relevant degenerations are {\it test configurations} which are defined as follows. For an integer $m>0$ a test configuration of exponent $m$ for $(X,L)$ is  a flat family of schemes $\pi:{\cal X}\rightarrow \bC$, with a relatively ample line bundle ${\cal L}\rightarrow {\cal X}$ and a $\bC^{*}$-action on ${\cal X},{\cal L}$ covering the standard action on $\bC$. For $t\in \bC$ we write $X_{t}$ for the scheme-theoretic fibre $\pi^{-1}(t)$ and $L_{t}$ for the restriction of ${\cal L}$ to $X_{t}$. We require that for all non-zero $t$ the pair $(X_{t}, L_{t})$ is isomorphic to $(X,L^{m})$. (Note that, due to the $\bC^{*}$-action, it suffices to know this for {\it some} non-zero $t$.) For technical reasons we also suppose that the total space ${\cal X}$ is normal.

The numerical criterion, in our situation, is provided by the Futaki invariant. In its general form, \cite{kn:D0}, this is defined for any $n$-dimensional projective scheme $Z$ with $\bC^{*}$-action and $\bC^{*}$-equivariant ample line bundle $\Lambda\rightarrow Z$, as follows. For each integer $k\geq 0$ we have a vector space  $H^{0}(Z;\Lambda^{k})$ with an induced $\bC^{*}$ action. Write $d(k)$ for the dimension of the space and $w(k)$ for the sum of the weights of the action. For large $k$,  $d(k)$ is given by a Hilbert polynomial which has degree exactly $n$ (since $\Lambda$ is ample), while $w(k)$ is given by a polynomial of degree at most $n+1$. Thus  $d(k)/kw(k)$  has an expansion, for large $k$:
$$   \frac{d(k)}{k w(k)} = F_{0} + k^{-1} F_{1} + \dots, $$
and the Futaki invariant is defined to be the co-efficient $F_{1}$. In our situation, we define the Futaki invariant ${\rm Fut}({\cal X})$ of a test configuration of exponent $m$ to be $m^{-1}$ times the Futaki invariant of the central fibre, with line bundle $L_{0}$ and the induced $\bC^{*}$-action. 

With these definitions in place we can state the main definition of this section.
\begin{defn}
  A polarised variety $(X,L)$ is K-semistable if  for any test configuration ${\cal X}$ we have ${\rm Fut}({\cal X})\geq 0$. It is $K$-stable if equality holds only when ${\cal X}$ is a product $X\times \bC$.
\end{defn}  
Note that in the last clause we allow the $\bC^{*}$-action on $X\times \bC$ to be induced from a non-trivial action on $X$. What we have called K-stability is often called K-polystability in the literature.   The precise statement of the result mentioned in the previous section, verifying Yau's conjecture, is
\begin{thm} A Fano manifold $X$ admits a Kahler-Einstein metric if and only if  $(X,K_{X}^{-1})$ is K-stable.
\end{thm} 
Here the \lq\lq only if''is usually regarded as the easier direction and is due to Berman \cite{kn:Berm}, following related results of various authors The uniqueness of the metric, modulo holomorphic automorphisms,  is a relatively old result of Bando and Mabuchi \cite{kn:BM}. We will not say more about these results here but focus on the \lq\lq if'' direction.

To give some background to the technical aspects of the proofs sketched in Section 4 below we will now try to explain why Theorem 1 is plausible. First we go back to the definition of the Futaki invariant of $(Z,\Lambda)$ in the case when $Z$ is a manifold, which was in fact the original context for Futaki's definition \cite{kn:Fut}. Choose a K\"ahler metric $\omega$ on $Z$ in the class $c_{1}(\Lambda)$ preserved by the action of $S^{1}\subset \bC^{*}$. Viewing $\omega$ as  a symplectic structure, this action is generated by a Hamiltonian function $H$ on $Z$. Then the Futaki invariant can be given by a differential geometric formula
\begin{equation}     \int_{Z} (R-\hat{R}) H \frac{\omega^{n}}{n!}, \end{equation}
where $R$ is the scalar curvature of $\omega$ and $\hat{R}$ is the average value of $R$ over $Z$. This formula can be derived from the equivariant Riemann-Roch theorem and can also be understood in terms of the asymptotic geometry of sections of $L^{k}$ as $k\rightarrow \infty$, in the vein of quasi-classical asymptotics in quantisation theory. What this formula shows immediately is that if $\omega$ can be chosen to have constant scalar curvature---in particular if it is a K\"ahler-Einstein metric---then the Futaki invariant vanishes. This given another way, different from the Matshusima theorem,  of ruling out K\"ahler-Einstein metrics on 1 or 2 point blow-ups of $\bC\bP^{2}$. The  definition of K-stability employs the Futaki invariant in a more subtle way; it is not just the  automorphisms of $X$ which need to be considered but of the degenerations.  The {\it Mabuchi functional} gives a way to understand this phenomenon.  This is a functional ${\cal F}$ on the space ${\cal H}$ of K\"ahler metrics  in a given cohomology class on a manifold $X$ defined via its first variation
\begin{equation}   \delta {\cal F} = \int_{X} (R-\hat{R}) \delta \phi \frac{\omega_{\phi}^{n}}{n!}. \end{equation}
Here $\delta \phi$ is an infinitesimal variation in the K\"ahler potential and one shows that such a functional ${\cal F}$ is well-defined, up to the addition of an arbitrary constant.  By construction a critical point of ${\cal F}$ is exactly a constant scalar curvature metrics which, in the setting of Theorem 1 can be shown to be K\"ahler-Einstein. (We  mention here that there is another functional, the {\it Ding functional} which has many similar properties to the Mabuchi functional and plays an important part in many developments. This is only defined for manifolds polarised by $K^{\pm 1}$.) 

 There are three possibilities:
\begin{itemize}
\item ${\cal F}$ is bounded below on ${\cal H}$ and attains its infimum;
\item ${\cal F}$ is bounded below but does not attain its infimum;
\item ${\cal F}$ is not bounded below.
\end{itemize}
An extension of Theorem 1 is the statement that these three possibilities correspond to $X$ being respectively $K$-stable, $K$-semistable (but not $K$-stable) and not $K$-semistable.

 Now suppose that ${\cal X}$ is a test configuration for $(X,K_{K}^{-1})$ and, for simplicity,  that the total space is smooth. Choose a K\"ahler metric on this total space, invariant under $S^{1}\subset \bC^{*}$. Pulling back by the $\bC^{*}$-action the restrictions of this metric to the fibres $X_{t}$ for non-zero $t$ can be regarded as a family of metrics $\omega(t)$ on the fixed manifold $X$ parametrised by $t\in \bC^{*} $ but these metrics have no limit, among metrics on $X$, as $t\rightarrow 0$. It is natural to think of this limit as a \lq\lq point at infinity'' in the space ${\cal H}$ of K\"ahler metrics on $X$. As we discuss further in Section 4.4 below,  the role of the Futaki invariant is to determine the asymptotic behaviour as $t\rightarrow 0$ of the Mabuchi functional in such families obtained from test configurations. 
 Theorem 1 can be understood roughly as saying that if there is no minumum of ${\cal F}$ this can be detected by studying the asymptotics at points at infinity of this kind (derived from algebro-geometric data).

Apart from its intrinsic interest---in giving an algebro-geometric criterion for the solubility of a PDE---Theorem 1 result also has implications for the construction of compactified moduli spaces of Fano manifolds, see work of Odaka \cite{kn:O2},  Spotti, Sun and Yao \cite{kn:SSY} and Li, Wang and Xu \cite{kn:LWX0}, \cite{kn:LWX}. For moduli questions, the notion of K-stability  is also relevant in the negative case of varieties with ample canonical bundle. In this case Odaka \cite{kn:O1} showed that a variety is K-stable if and only if has semi log canonical singularities, which is equivalent to stability in the sense of Alexeev, K\'ollar and Shepherd-Barron. Berman and Guenancia showed that for such varieties this is equivalent to the existence of a K\"ahler-Einstein metric \cite{kn:BG} (with a suitable definition in the singular case). Compactifying moduli spaces of manifolds involves adding points corresponding to singular varieties and it is interesting to relate the behaviour of the K\"ahler-Einstein metrics to the algebraic geometry of the singularities.  There is much recent progress in this direction, see \cite{kn:HS} for example.

\section{Riemannian convergence theory and projective embeddings}

In this section we will discuss some ideas which play an important role in three of the four proofs considered in Section 4 below. The general context can be explained as follows. In solving a PDE problem {\it compactness}---the ability to take limits in some kind of approximating scheme---is usually crucial. On the other hand in our problem we need to exhibit the obstruction to solving the problem (the existence of a K\"ahler-Einstein metric) as an algebro-geometric object (a test configuration with non-positive Futaki invariant). In the framework of Ricci curvature in Riemannian geometry there is a well-developed convergence theory of {\it Gromov-Hausdorff limits}; thus, in the K\"ahler situation, we would like to relate these limits to algebraic geometry and that is the topic of this section. 
  
  We begin by recalling some of the main results from the Riemannian theory of manifolds with a lower bound on the Ricci curvature. The foundation of the theory is the link between the Ricci curvature and volume expressed by the Bishop comparison theorem. For simplicity we just consider the case of an $m$-dimensional Riemannian manifold $M$ with ${\rm Ricci}\geq 0$. Then Bishop's theorem states that for each $p\in M$  the volume ratio
\begin{equation}  v_{p}(r) = \frac {{\rm Vol}(B_{p}(r))}{\Omega_{m} r^{m}},\end{equation}
is a weakly decreasing function of $r$. Here $B_{p}(r)\subset M$ is the metric ball of radius $r$, and we introduce the normalising constant $\Omega_{m}$---the volume of the unit ball in $\bR^{m}$---so that $v_{p}(r)$ tends to $1$ as $r$ tends to $0$. If $M$ is  compact with total volume $V$ and diameter $\leq D$ it follows that
\begin{equation}   {\rm Vol}(B_{p}(r)) \geq  \kappa r^{m}, \end{equation}
with $\kappa= V/D^{m}$. Recall that the {\it Gromov-Hausdorff distance} between two compact metric spaces $A,B$ is defined as the infimum of the numbers $\delta$ such that there is a metric on the disjoint union $A\sqcup B$ which extends the given metrics on $A,B$ and such that both $A,B$ are $\delta$-dense in $A\sqcup B$. If $(M_{i}, g_{i})$ is a sequence of compact Riemannian $m$-manifolds with ${\rm Ricci}\geq 0$, ${\rm Vol}(M_{i})=V$, ${\rm diam}(M_{i})\leq D$ then Gromov's compactness theorem asserts that there is a subsequence which converges in the sense of this Gromov-Hausdorff distance to some limiting metric space $(Z, d_{Z})$. (More generally, the same result applies if we have any fixed lower bound on the Ricci curvatures.) The proof is an elementary argument based on ball-packing considerations and the lower bound (6). It is sometimes convenient to express  this Gromov-Hausdorff convergence in terms of a  natural topology on the disjoint union
 $$        {\cal M}= Z\cup \bigsqcup_{i}M_{i} . $$
 Thus for $q\in Z$ it makes sense to talk about points $p_{i}\in M_{i}$ which are close to $q$. 

Results of Anderson \cite{kn:And}, Cheeger-Colding \cite{kn:CC} and Cheeger-Colding-Tian\cite{kn:CCT} give finer information about such \lq\lq non-collapsed'' Gromov-Hausdorff limits. (Here the non-collapsing refers to the volume lower bound, which rules out the collapse of the sequence of $m$-dimensional manifolds to some lower dimensional space.) The notion of Gromov-Hausdorff convergence can be extended to sequences of spaces with base points: the metric balls of any fixed radius centred at the base points are required to converge as above. For each point $q\in Z$ and sequence of real numbers $\lambda_{j}\rightarrow \infty$ we consider the sequence of based metric spaces $(q, Z, \lambda_{i} d_{Z})$. After perhaps passing to a subsequence  we have a based Gromov-Hausdorff limit which is a metric cone, a {\it tangent cone} of $Z$ at $q$. The regular set $R\subset Z$ is defined to be the subset where some tangent cone is ${\bf R}^{m}$ and the complement $Z\setminus R$ is the singular set $S$. If the manifolds $M_{i}$ satisfy a fixed bound on the Ricci curvature $\vert {\rm Ricci}\vert \leq \Lambda$ then more can be said. Anderson showed that there are fixed $\delta_{m},\kappa_{m}$, depending only the dimension, such that if $p\in M_{i}$ and $r\leq \Lambda^{-1/2}$ then if the if the volume ratio $v_{p}(r)$ is greater than  $1-\delta_{m}$ there are harmonic co-ordinates on the sub-ball  $B_{p}(\kappa r)$ in which the metric satisfies $C^{1,\nu}$ estimates. Together with results of Cheeger and Colding this shows that the regular set is open in $Z$ and the Riemannian metrics converge to a  $C^{1,\nu}$-Riemannian metric $g_{\infty}$ on $R$. ( This means that  if $U$ is a pre-compact open set in $R$ there are $C^{2,\nu}$ diffeomorphisms $\chi_{i}:U\rightarrow M_{i}$ which, regarded as maps into ${\cal M}$,  converge to the inclusion $U\rightarrow Z$ and such that the $\chi^{*}_{i}(g_{i})$ converge in $C^{1,\nu}$ to $g_{\infty}$.) The singular set $S$ has Hausdorff codimension at least $4$. (In the general Riemannian context, this is a recent result of Cheeger and Naber \cite{kn:CN}, but in the K\"ahler case which will be our concern it goes back to Cheeger, Colding and Tian \cite{kn:CCT}.) The corresponding statements apply to tangent cones: each has a smooth Ricci-flat metric outside a closed singular set of codimension  at least 4. 

We want to relate these ideas to algebraic geometry and in this section we will focus on the case considered in \cite{kn:DS}. Thus we suppose that $(X_{i}, \omega_{i})$ are K\"ahler manifolds of complex dimension $n$ with fixed volumes $V$, diameters $\leq D$ and $\vert {\rm Ricci}\vert\leq \Lambda$ and with a Gromov-Hausdorff limit $Z$, as above. We suppose that that these are polarised manifolds, so that $\omega_{i}$ is the curvature of a Hermitian holomorphic line bundle $L_{i}\rightarrow X_{i}$. The main result is that $Z$ can be endowed with the structure of an algebraic variety. More precisely,  we allow passage to a subsequence and form the disjoint union ${\cal M}$ as above. Then there is a continuous map $I:{\cal M}\rightarrow \bC\bP^{N}$ with the two properties.
\begin{itemize}\item There is some fixed $k$ such that for sufficiently large $i$, the restriction of $I$ to $X_{i}$ is an embedding defined by the holomorphic sections of $L_{i}^{k}\rightarrow X_{i}$;
\item The restriction of $I$ to $Z$ is a homeomorphism to its image, which is a normal projective variety in $\bC\bP^{N}$. 
\end{itemize}

This result can be seen as an extension of the Kodaira embedding theorem to   singular limit spaces and the proof extends some the ideas in one  approach to the Kodaira theorem (an approach which seems to be well-known to experts but does not feature in standard textbooks). Suppose that $L\rightarrow X$ is a holomorphic line bundle over a compact complex manifold and $\sigma_{0}$ is a holomorphic section of $L$ over an open subset $V\subset X$. Let $\beta$ be a cut-off function with compact support in $V$, extended by $0$ over $X$ and with $\beta=1$ on some interior region $V_{0}\subset V$. Then we can regard $\sigma=\beta \sigma_{0}$ as a smooth section of $L$ over $X$ in an obvious way. This will not be  a holomorphic section because we have a term coming from the cut-off function
$$    \db \sigma = (\db \beta) \sigma_{0}, $$
 but if we have  Hermitian metrics on $L$ and $X$ we can project $\sigma$ to the space of holomorphic sections using the $L^{2}$ inner product on sections of $L$, arriving at a holomorphic section $s$. Of course in this generality the construction need not be useful---the projection $s$ could be $0$. The idea is that under suitable hypotheses we can arrange that $s$ is not zero and is very close to the original section $\sigma_{0}$ over $V_{0}$. The key  is to estimate the error term $\eta=\sigma-s$, which is given by a Hodge Theory formula
$$   \eta= \db^{*} \Delta^{-1}  \db \sigma, $$ where $\Delta$ is the $\db$-Laplacian on the $L$-valued $(0,1)$ forms. Of course this formula only makes sense if $\Delta$ is invertible, i.e. if the cohomology group $H^{1}(X;L)$ is zero. Then we have $$ \Vert \eta\Vert^{2}_{L^{2}}= \langle \db^{*} \Delta^{-1} \db \sigma, \db^{*}\Delta^{-1} \db \sigma\rangle= \langle \db\db^{*}\Delta^{-1}\db \sigma, \Delta^{-1}\db\sigma\rangle=\langle \db\sigma, \Delta^{-1}\db\sigma\rangle, $$ and so
$$  \Vert \eta\Vert_{L^{2}}^{2}\leq \Vert \Delta^{-1}\Vert \ \Vert \db\sigma\Vert^{2}, $$ where $\Vert \Delta^{-1}\Vert$ is the $L^{2}$-operator norm.  Now suppose that $L$ is a positive line bundle and the metric on $X$ is the K\"ahler metric $\omega$ given by the curvature of $L$. There is then a formula of Weitzenb\"ock type
    \begin{equation}   \Delta = P^{*} P + {\rm Ric} + 1 , \end{equation}
    where $P$ is the $(0,1)$-component of the covariant derivative on $L$-valued $(0,1)$ forms. But all we use is that $P^{*}P$ is a non-negative operator, so if ${\rm Ric} \geq -1/2$ (say) then $\Delta\geq 1/2$ and $\Delta^{-1}$ is defined and has $L^{2}$-operator norm at most $2$. So we get $$\Vert \eta\Vert_{L^{2}}\leq \sqrt{2} \Vert \db \sigma \Vert_{L^{2}}.$$
 So if we can arrange that $\db\sigma$ is small, in $L^{2}$ norm compared with $\sigma$,  then we get a non-trivial holomorphic section $s$. By construction, $\eta$ is holomorphic over $V_{0}$ and the $L^{2}$ norm of $\eta$ controls all derivatives there, by the usual elliptic estimates,  and we can hope to show that $s$ is a small perturbation of $\sigma_{0}$ in any $C^{r}$ norm. The whole approach extends easily to a case when the original section $\sigma_{0}$ is not exactly holomorphic but approximately so, measured in terms of bounds on the $L^{2}$  norm of $\db \sigma_{0}$. 

To prove the usual Kodaira embedding theorem in this way we first make the simple observation that changing $L$ to $L^{k}$ for some $k>1$ corresponds to rescaling the metric by a factor $k$---i.e. scaling lengths of paths by a factor $\sqrt{k}$. Under this rescaling the Ricci curvature is multiplied by $k^{-1}$ so if we start with any positive line bundle and take a suitable power we can arrange that the condition ${\rm Ric}\geq - 1/2$ holds after rescaling.  We take $U\subset X$ to be a small co-ordinate ball centred on a point $p\in X$. The flat model is given by the trivial line bundle $\Lambda$ over $\bC^{n}$ with a Hermitian structure corresponding to the Euclidean metric on $\bC^{n}$. In  this structure the trivialising section $\tau$ of $\Lambda$ has Gaussian norm 
   $$  \vert \tau(z) \vert = e^{-\vert z\vert^{2}/4}, $$
   which decays very rapidly at infinity. After rescaling,  the geometry of the manifold $X$ in the small ball $U$ is close to the flat model and we get an approximately holomorphic section $\sigma_{0}$ modelled on $\tau$. By suitable choice of parameters one arranges a cut-off function $\beta$ with the support of  $\nabla \beta$ contained in the region where $\sigma_{0}$ is very small,  so that $\Vert \db \sigma\Vert$ is also very small. Here the \lq\lq suitable choice of parameters'' involves choosing $k$ large. The upshot is that for a suitable $k$ and for each $p\in X$ we construct a  holomorphic section of $L^{k}$ \lq\lq peaked'' around $p$ and in particular not vanishing at $p$. Thus the sections of $L^{k}$ define a holomorphic map of $X$ to a projective space and one can go further to show that this map is an embedding, once $k$ is sufficiently large.  (Note that the formula (7) is the same as that used in the proof of the Kodaira-Nakano vanishing theorem, which is invoked in the usual proof of  the embedding theorem via blow-ups.)

  Returning to our main discussion, we do not want to do analysis directly on the Gromov-Hausdorff limit $Z$ but instead establish uniform estimates on the converging sequence $X_{i}$. Consider any polarised manifold $(X,L)$ with K\"ahler metric given by the curvature of $L$. We endow $H^{0}(X,L)$ with the standard $L^{2}$ metric and for $x$ in $X$ we consider the evaluation map
$${\rm ev}_{x}: H^{0}(X,L)\rightarrow L_{x}.$$
We define $\rho_{L}(x)$ to be square of the norm of this map so the statement that $\rho_{L}(x)>0$ for all $x\in X$ is the same as saying that the sections of $L$ define a map  $\tau: X\rightarrow \bP(H^{0}(X,L)^{*})=\bP$. More generally, a lower bound on the $\rho_{L}(x)$ gives metric control of this map. The operator norm of
   $$  d\tau_{x}: TX_{x}\rightarrow T\bP_{\tau(x)}$$ is
 $\rho_{L}(x)^{-1/2} \max (\vert (\nabla s)_{x}\vert)  $, where the maximum is taken over holomorphic sections $s$ of $L$ with $L^{2}$ norm $1$ vanishing at $x$. In our situation,  with Ricci curvature and diameter bounds, there is a well-known upper bound on $\vert \nabla s\vert$, so a strictly positive lower bound on the $\rho_{L}(x)$ gives a Lipschitz bound on the map $\tau$. Replacing $L$ by $L^{k}$ we see that the crucial point is to  find some $k$ and $b>0$ so that for all $i$ and for all $x$ in $X_{i}$ we have 
\begin{equation}  \rho_{L^{k}}(x)\geq b. \end{equation}
(Such a bound is sometimes referred to as a {\it partial $C^{0}$-estimate}.) It is straightforward to reduce to the case when the dimension of $H^{0}(X_{i},L_{i}^{k})$ is independent of $i$, so that we can regard $\tau_{i}:X_{i}\rightarrow \bP$ as mapping into a fixed projective space. If the bound (8) holds then, after possibly taking a subsequence,  we can pass to the Gromov-Hausdorff limit and define a continuous map $\tau_{\infty}:Z\rightarrow \bP$ with image a projective variety. Further arguments then show that after perhaps increasing $k$ this map $\tau_{\infty}$ is an homeomorphism from $Z$ to a normal projective variety. 

The central issue then is to establish the lower bound (8), as has been emphasised in many places by Tian. Let $q$ be a point in $Z$ and $C(Y)$ be some tangent cone to $Z$ at $q$. If $q$ is a smooth point then $C(Y)$ is $\bC^{n}$ and we have the model Gaussian section as discussed above. In general we always have a Hermitian line bundle $\Lambda$ over the regular part of $C(Y)$ with a holomorphic section $\tau$ satisfying exactly the same Gaussian decay  with respect to the distance to the vertex of the cone. We choose a suitable open subset $U$ of the regular part of the cone. It follows from the definitions that for large $i$ there are diffeomorphisms $\chi_{i}:U\rightarrow X_{i}$ which are approximately holomorphic isometries with respect to rescalings of the metrics $k\omega_{i}$ for some suitable large $k$. (Here the approximation can be made as close as we like by taking $k$ large.) Then are then two main technical points to address.
\begin{itemize}
    \item We want to have lifts $\tilde{\chi}_{i}:\Lambda \rightarrow \chi^{*}_{i}(L^{k})$ to approximate isomorphisms of Hermitian line bundles over $U$. 
\item We want a suitable cut-off function $\beta$ on $U$ with $\vert \db \beta  \tau\vert$ small in $L^{2}$. 
\end{itemize}
Given these, we can transport the section $\beta \tau$ to an approximately holomorphic section of $L_{i}^{k}$ over $\chi_{i}(U)$and follow the projection procedure to get a holomorphic section $s$ of $L_{i}$ modelled on $\tau$. The derivative bounds on $s$ give a lower bound on $\vert s\vert$ over all points in $X_{i}$ close to $q$ and an elementary covering argument establishes the bound (8).

The first technical point involves considerations of the holonomy of the connection on $\Lambda^{*}\otimes \chi_{i}^{*}(L_{i}^{k})$, which has very small curvature by construction--- this is straightforward if $U$ is simply connected. The second technical point involves the singular set in $C(Y)$, and in particular the fact that this has Hausdorff codimension strictly greater than $2$ (see 4.1 below).

\section{Four proofs}
We now come to the core of this survey in which we discuss four different proofs of the equivalence between stability and K\"ahler-Einstein metrics. In  total these proofs run to many hundreds of pages so it is  impossible  to give any kind of thorough account of them here. All we can do is to explain general strategies and some salient points in the arguments.

\subsection{ The proof by cone singularities}

(Note that the announcement \cite{kn:CDS0} contains an outline of this proof.)  Given a Fano manifold $X$ we fix some suitable $m\geq 1$ and a smooth divisor $D$ in $\vert -m K_{X}\vert$. For $0<\beta\leq 1$ we can define a class of  K\"ahler metrics on $X$ with cone singularity of angle $2\pi \beta$ along $D$ and extend the whole theory to this case. (When $\beta= r^{-1}$ for an integer $r$ these are  orbifold metrics, and hence well-established.  There are also close analogies with the theory of parabolic structures and singular Hermitian  Yang-Mills connections as developed in \cite{kn:Biq} for example.) We can define a modified Mabuchi functional 
\begin{equation}\delta {\cal F}_{\beta}= \int_{X} (R-\hat{R}) \delta \phi + (1-\beta) \int_{D} \delta \phi - c \int_{X}\delta \phi \end{equation}
where $c= (m c_{1}(X)^{n})^{-1}$ is the ratio of the volume of D to the volume of X, so that the right hand side vanishes when $\delta \phi$ is a  constant.  Roughly speaking, we have a family of functionals ${\cal F}_{\beta}$ with critical points the K\"ahler-Einstein metrics with cone angle $\beta$ and the strategy of the proof is to follow a family of such critical points as $\beta$ increases. We want to show that either the family continues up to $\beta=1$, which gives our desired K\"ahler-Einstein metric on $X$,  or that the critical point moves off to infinity and that this yields a test configuration violating the K-stability condition.  

To begin we need to show that a solution exists for {\it some} $\beta$. Take $m>1$ and $\beta= r^{-1}$ so that we are in the orbifold case. If $r>m/m-1$ then characteristic class arguments show that the we are in the situation of negative Ricci curvature and the desired solution follows from a straightforward orbifold extension of the standard Aubin-Yau theory. Next one shows that the set of $\beta$ for which a solution exists is {\it open}. This can be achieved using a suitable linear elliptic theory on manifolds with cone singularities \cite{kn:D1}. 

Suppose then that there is some $\beta_{0}\in (0,1]$ such that a solution $\omega_{\beta}$ exists for $\beta<\beta_{0}$ but that there is no solution for $\beta=\beta_{0}$. We need to extend the theory sketched in Section 3, for smooth metrics with bounded Ricci curvature, to  K\"ahler-Einstein metrics with cone singularities. To begin we show that a metric with cone singularity can be approximated in the Gromov-Hausdorff  sense by smooth metrics with Ricci curvature bounded below (\cite{kn:CDS}, Part I). Then the Cheeger-Colding theory implies that there is a subsequence $\beta_{i}$ increasing to $\beta_{\infty}$ and a Gromov-Hausdorff limit $Z$ of the $(X,\omega_{\beta_{i}})$ and $Z$ has metric tangent cones at each point. We call a tangent cone $C(Y)$ \lq\lq good''  if the regular set $ C(Y_{\reg})$ is open, the metric is induced by a smooth K\"ahler metric  there and  and for each compact subset $K\subset Y_{ \reg}$ and each $\eta>0$ there is a cut-off function $\gamma$ of compact support in $Y_{\reg}$, equal to $1$ on $K$ and with
\begin{equation}  \Vert \nabla \gamma\Vert_{L^{2}}\leq \eta. \end{equation}
The main technical result is that all tangent cones to $Z$ are good. Given this, an extension of the arguments outlined in Section 3 above show that $Z$ is naturally a normal projective    variety, carrying a singular K\"ahler-Einstein metric $\omega_{\infty}$. Moreover if we write $X_{i}$ for the metric space $(X,\omega_{\beta_{i}})$ and $D_{i}\subset X_{i}$ for the divisor $D$ then there is a divisor $\Delta\subset Z$ such that the pairs $(X_{i}, D_{i})$ converge to $(Z,\Delta)$.  

The new feature in this case, which leads to the difficulty in proving that tangent cones are \lq\lq good'', is the possibility of codimension 2 singular sets. This is the critical dimension with respect to the cut-off control (10). If $\psi$ is a compactly-supported function on $\bR^{m}$ and if $\psi_{\lambda}$ is the rescaled function $\psi_{\lambda}(x)=\psi(\lambda^{-1} x)$ then
$$   \Vert \nabla \psi_{\lambda}\Vert_{L^{2}} = O(\lambda^{(m-2)/2})$$
which tends to $0$ with $\lambda$ if $m>2$. This allows one to construct cut-off functions with derivative arbitrarily small in $L^{2}$ adapted to a compact set $A$ of Hausdorff codimension strictly greater than $2$. In the codimension 2 case one needs appropriate control of the volume of the $\lambda$-neighbourhood $N_{\lambda}(A)$:
$$ {\rm Vol}\  N_{\lambda}(A))\leq C \lambda^{2}. $$
This is equivalent to the notion of {\it Minkoswki codimension}$\geq 2$ i.e.
for any $r$ the set $A$ can be covered by $O(r^{-2})$ balls of radius $r$.

To complete the proof of the Theorem one wants to show that if indeed the family of solutions breaks down at some $\beta_{\infty}$ as considered above then there is test configuration with central fibre $Z$ and with non-positive Futaki invariant, so the original manifold $X$ is not K-stable. To this end one can extend the whole theory of stability and test configurations to pairs consisting of a variety and divisor, with a real parameter $\beta$.  There is a modified Futaki invariant ${\rm Fut}_{\beta}$  which compares with the usual formula (3) (in the smooth case) just as the modified Mabuchi functional (9) compares with (4). One wants to construct a test configuration $({\cal X}, {\cal D})$ for the pair $(X,D)$ with central fibre $(Z,\Delta)$ and show that
\begin{equation}   {\rm Fut}_{\beta_{\infty}}({\cal X}, {\cal D})= 0. \end{equation}
The Futaki invariant ${\rm Fut}_{\beta}$ depends linearly on $\beta$ so the fact that $(X,D)$ is stable for small $\beta$ implies that ${\rm Fut}({\cal X})={\rm Fut}_{1}({\cal X}, {\cal D})\leq 0$. 

Let $G$ be the automorphism group of the pair $(Z,\Delta)$--a complex Lie group. The existence of this test configuration $({\cal X}, {\cal D})$ follows from general principles once it is established that $G$ is {\it reductive}; the complexification of a compact subgroup $K\subset G$. Using  projective embeddings in some $\bC\bP^{N}$,  the pairs  correspond to points $[X,D]$, $[Z,\Delta]$ in a suitable Hilbert scheme ${\bf S}$ which in turn is embedded in some large projective space $\bP$.  The group $PGL(N+1,\bC)$ acts on ${\bf S}$ and $\bP$ and what we know from the convergence discussion above is that $[Z,\Delta]$ is in the closure of the orbit of $[X,D]$. The group $G$ can be identified with the stabiliser of $[Z,\Delta]$ in $PGL(N+1,\bC)$. If $G$ is reductive a version of the Luna slice theorem gives a slice for the action of $PGL(N+1,\bC)$ at $[Z,\Delta]$ (one takes a $G$-invariant complement to the action of $G$ on the tangent space of $\bP$ at $[Z,\Delta]$). A well-known result of Hilbert and Mumford, applied to the $G$ action on this slice, shows that there is a $1$-parameter subgroup  $\bC^{*}\subset G$ such that $[Z,\Delta]$ is in the closure of the $\bC^{*}$-orbit of $[X,D]$ and this is equivalent to the desired test configuration. 

The reductivity of $G$ is an extension of the Matsushima result to the case of pairs and singular varieties.  Likewise the vanishing of the Futaki invariant (11) is an extension of the simple case  considered in Section 2 above, for manifolds of constant scalar curvature. For the proofs one has to work with the singular metric $\omega_{\infty}$ using techniques from pluripotential theory.

\subsection{The proof by the continuity method}

In the continuity method one fixes a positive $(1,1)$ form $\alpha$ representing $c_{1}(X)$ and tries to solve the family of equations for $\omega_{s}$, with parameter $s\in [0,1]$:

\begin{equation}   {\rm Ricci}(\omega_{s})= (1-s) \alpha + s \omega_{s}. \end{equation}

Yau's solution of the Calabi conjecture shows that there is a solution for $s=0$ and it is well-known that the set of parameters for which a solution exists is open; if the solution can be continued up to $s=1$ we have a K\"ahler Einstien metric, so the problem is to prove closedness (under the stability hypothesis). Thus we suppose that for some $S<1$ there are solutions for $s<S$ but none for $s=S$. The set-up is very similar to that of cone singularities above, indeed the latter can be regarded as a variant of the continuity method, replacing the smooth form $\alpha$ with the current defined by a divisor.

If one knew that the fixed form $\alpha$  was bounded with respect to $\omega_{s}$ then the $\omega_{s}$ would have bounded Ricci curvature and the results discussed in Section 3 above would apply immediately to give a limiting metric on a normal projective variety. So the major difficulty is that we do not have such a bound. The Ricci curvature of the $\omega_{s}$ is positive so the fundamentals of the Cheeger-Colding theory apply and we obtain a Gromov-Hausdorff limit $Z$ (more precisely,  a limit of some sequence $(X, \omega_{s_{i}})$ with $s_{i}$ increasing to $S$). But this theory does not ensure that the regular set in $Z$ is open, or give the good convergence properties over the regular set exploited in the argument of Section 3. This is one of the main problems overcome by Sz\'ekelyhidi in \cite{kn:Gabor2}.

  To explain some of Sz\'ekelyhidi's arguments, we restrict attention to the case when $S<1$.   Consider first a unit ball $B$ centred at a point $p$ in a K\"ahler manifold,  with metric $\omega$,  and a vector-valued holomorphic function
$  f:B\rightarrow \bC^{m}$. Suppose that the Ricci curvature is bounded in $B$, say $\vert {\rm Ric}\vert \leq 4$, and that the pair $(\omega, f)$ satisfy the equation
\begin{equation}    {\rm Ric}(\omega) = \sigma \omega + f^{*}(\OmegaEuc), \end{equation}
for some $\sigma\in [0,1]$ where $\OmegaEuc$ is the standard Euclidean K\"ahler form on $\bC^{m}$. We claim  that for any $\epsilon>0$ we can find a $\delta$ (independent of $\sigma$) such that for any such $\omega$,  if the volume ratio $v_{p}(1)$  exceeds $(1-\delta)$ then in fact
$$\vert {\rm Ric}(p)\vert \leq \epsilon.   $$ 

First, the Ricci curvature is non-negative so by the Bishop inequality we can pass to a smaller ball with centre $p$ (rescaled) and preserve the volume bound. By the results of Anderson we may as well assume that the metric on $B$ is $C^{1,\nu}$-close to the Euclidean metric in harmonic co-ordinates. By a suitable version of the Newlander-Nirenberg integrability theorem we can also suppose that these co-ordinates are actually holomorphic. The bound on the Ricci curvature and the equation (13) means that $\vert \nabla f\vert^{2} \leq 2n $ and since $f$ is holomorphic we get interior bounds on all higher derivatives of $f$ and hence on the Ricci tensor,  in these holomorphic co-ordinates. Thus if $\vert{\rm Ricci}(p)\vert>\epsilon$ we will have $\vert {\rm Ricci}\vert>\epsilon/2$ over a ball of definite size centred at $p$. The $C^{1,\nu}$ bound on the metric tensor gives $C^{,\nu}$ bounds on the Christoffel symbols. From this Sz\'ekelyhidi shows that there is some unit tangent vector $v$ at $p$ such that, in geodesic polar-coordinates centred at $p$,  there is a definite lower bound on
    ${\rm Ric} (\frac{\partial}{\partial r}, \frac{\partial}{\partial r})$ at all points close to $p$ and along geodesics starting from $p$ at a sufficiently small angle to $v$. Then the proof of the Bishop inequality shows that this Ricci curvature {\it reduces} the volume of $B$ by a definite amount (compared with the Euclidean ball) determined by $\epsilon$. So ${\rm Vol}(B)\leq (1-\delta(\epsilon))\Omega_{2n}$, say. Choosing $\delta=\delta(\epsilon)$ we get a contradiction to the hypothesis that $\vert {\rm Ric}(p)\vert\geq \epsilon$ and the claim is established. 

Clearly the result extends (with a suitable $\delta(\epsilon)$) to the case when $\OmegaEuc$ is replaced by any smooth positive $(1,1)$-form $A$ defined over a suitable neighbourhood of the image of $f$. In the case at hand one can cover $X$ by a finite number of holomorphic co-ordinate charts. Working near a given point in $X$ such a chart yields the holomorphic map $f$ above (with $m=n$) and we take the form $A$ corresponding to $(1-s)\alpha$ in this chart. If $s\leq S<1$ we get a $\delta(\epsilon)$ such that the discussion above applies to any rescaling of a small ball in $(X,\omega_{s})$.

Now let $q$ be a point in the regular set of the Gromov-Hausorff limit $Z$ and let $p_{i}\in (X,\omega_{s_{i}})$ be a sequence converging to $q$ in the sense of the Gromov-Hausdorff convergence. Let $B_{i}$ be the unit ball obtained by rescaling a small ball (of fixed radius $\rho$, independent of $i$) about $p_{i}$. For any given $\delta>0$ we can suppose that for all subballs $\tilde{B}\subset B_{i}$ (not necessarily centred at $p_{i}$) the \lq\lq volume defect'' of $\tilde{B}$ is less than $\delta$. We choose $\delta$ as above, for some $\epsilon<1$. Now let
$$M={\rm max}_{x\in B_{i}} \left(\vert {\Ric}(x)\vert\  d(x,\partial B_{i})^{2}\right). $$ A standard line of argument shows that in fact $M\leq 4$. For if not 
 let $\tilde{p}\in B_{i}$ be a point where the maximum is attained and $\tilde{B}$ be the ball of radius $\tilde{d}/\sqrt{M}$ centred at $\tilde{p}$, where $\tilde{d}$ is  the distance from $\tilde{p}$ to the boundary of $B_{i}$. Rescaling $\tilde{B}$ to unit size we get a ball to which the previous results apply. After rescaling those results give $\vert {\rm Ricci}(\tilde{p})\leq \epsilon d_{0}^{-2} M< d_{0}^{-2}M$ which is a contradiction to the choice of $\tilde{x}$.  Thus $M\leq 4$  and in particular $\vert {\rm Ricci}(p_{i})\vert \leq 4 \rho^{-2}$. The conclusion  is that Sz\'ekelyhidi is able to show that, when $S<1$ the Ricci curvature is bounded near points in the regular set in $Z$. It follows that the regular set is open and  carries a $C^{1,\nu}$ K\"ahler metric. Going further, he extends the discussion to tangent cones and shows that these are all \lq\lq good'' in the sense discussed in 4.1 above.  Note that in this situation the singular set can have real codimension 2, different from the simpler situation considered in Section 3. (In fact results of C. Li \cite{kn:Li} in the toric case suggest that the limit as $s\rightarrow S$ will develop cone singularities along a divisor.) 

  In \cite{kn:Gabor2},  Sz\'ekelyhidi established the partial $C^{0}$ estimate along the continiuty method and used this to show that another notion of stability, introduced by S. Paul, implies the existence of a K\"ahler-Einstein metric. The corresponding result for $K$-stability was established in the subsequent paper \cite{kn:DaS} of Datar and Sz\'ekelyhidi. The output from the convergence theory is that if the continuity method breaks down at $S\in (0,1]$ there is a limiting projective variety $Z$, a singular K\"ahler metric $\omega_{S}$ and  a closed non-negative $(1,1)$ current  $\alpha_{\Psi}$ on $Z$ satisfying the equation
$$   {\Ric}(\omega_{S})= (1-S)\omega_{S}+ S \alpha_{\Psi}. $$
The possible presence of singularities means that that this equation needs to be interpreted. The current $\alpha_{\Psi}$ is locally written as $i\dbd\psi$ where $ \psi$ is an $L^{1}$ purisubharmonic function. Globally, $\Psi$ is a singular Hermitian metric on the anticanonical bundle. Datar and Sz\'ekelyhidi set up a theory for pairs consisting of a variety and a $(1,1)$-current, analogous to the theory for pairs (variety, divisor) discussed above. The new feature is that their space of pairs is infinite dimensional. They are able to carry through a similar strategy to that outlined in 4.1 above by approximating $(1,1)$-currents by those defined by divisors. The form $\alpha$ on $X$ can be taken to be the restriction of the Fubini-Study metric under an embedding $X\subset \bP$,  then there is a integral geometry formula
$$  \alpha= \int_{\bP^{*}}  [H\cap X] d\mu(H), $$
where $\bP^{*}$ is the dual projective space parametrising hyperplanes $H\subset \bP$, $\mu$ is the standard measure on $\bP^{*}$ and $[H\cap X]$ is the current of the divisor $H\cap X$ in $X$. Replacing the integral by a finite sum gives the approximation procedure which is the starting point for these arguments. 
 
The results of Datar and Sz\'ekelyhidi go further than the statement of Theorem 1 in two directions. First they prove an analogous result for solutions of the {\it K\"ahler-Ricci soliton} equation. Recall that this equation is
\begin{equation}  {\rm Ricci}(\omega) - \omega= L_{v}\omega, \end{equation}
where $v$ is a holomorphic vector field and $L_{v}$ is the Lie derivative. Such metrics are  the appropriate analogues of K\"ahler-Einstein metrics on Fano manifolds with non-vanishing Futaki invariant and represent fixed points of the K\"ahler-Ricci flow (modulo holomorphic diffeomorphisms), which we discuss further below.

In another direction, Datar and Sz\'ekelyhidi's proof is compatible with group actions so they prove that to test K-stability of a Fano manifold $X$ it suffices to consider test configurations with an additional compatible action of ${\rm Aut}(X)$. This is important because an outstanding defect of the general theory is that it is very hard to verify K-stability of a polarised variety. The problem becomes more tractable for manifolds with large symmetry groups. Toric manifolds, with a complex torus action having an dense orbit, can be described in terms of polytopes. Both K\"ahler metrics invariant under the action of the corresponding real torus  and toric test configurations  can be described by convex functions on the polytope and the stability condition is relatively explicit.  However in the toric Fano case the existence problem for K\"ahler-Einstein metrics and K\"ahler-Ricci solitons was completely settled by Wang and Zhu in 2004 \cite{kn:WZ}, and no interesting phenomena arise, from the point of view of stability. Ilten and S\"uss \cite{kn:IS} consider $n$-dimensional varieties with an action of an $(n-1)$-dimensional complex torus and develop a combinatorial description of these. In this way they are able to produce new examples, of manifolds which are K-stable, and the theorem of Datar and Szekelyhidi gives corresponding explicit new results about the existence of K\"ahler-Einstein metrics. In a similar vein, Delcroix studied group compactifications \cite{kn:Delcroix1}; that is he considered a manifold  $X$ which contains a complex reductive Lie group $G$ as a dense open subset and such that both left and right translations on $G$ extend to $X$. These can be described by polytopes in the Lie algebra of a maximal compact real torus in $G$ and Delcroix extends the arguments of Wang and Zhu to find an explicit condition for the existence of a K\"ahler-Einstein metric. This work of Delcroix was essentially self-contained and did not invoke general existence results such as Theorem 1. A subsequent paper \cite{kn:Delcroix2} later showed that his condition emerges from an analysis of equivariant degenerations of $X$ and extended the results to the larger class of spherical varieties, using the theorem of Datar and Sz\'ekelyhidi.

\subsection {The proof by Ricci flow}

If $X$ is a Fano mainfold the relevant version of the Ricci flow is the evolution equation
\begin{equation} \frac{\partial \omega_{t}}{\partial t} = \omega_{t} -{\rm Ricci}(\omega_{t}), \end{equation} for a one-parameter family of metrics $\omega_{t}$ in the class $c_{1}(X)$.
This can be expressed in terms of the K\"ahler potential. For each $t$ there is a unique $h_{t}$ such that
   $$   \omega_{t}-{\rm Ricci}(\omega_{t}) = i\dbd h_{t}, $$
   normalised so that $\max_{X} h_{t}=0$ and we can write $\omega_{t}=\omega_{0}+ i \dbd \phi_{t}$ where $\phi_{t}$ evolves by:
\begin{equation}    \frac{\partial\phi_{t}}{\partial t} = h_{t}. \end{equation}
 It has been known for many years that this equation has a solution for all $t\in[0,\infty)$, starting with any initial condition. The main result of the paper of Chen, Sun and Wang \cite{kn:CSW} is that this flow converges as $t\rightarrow \infty$ to a \lq\lq weak'' Kahler-Ricci soliton metric $\omega_{\infty}$ on a normal projective variety $X_{\infty}$ (in fact a \lq\lq Q-Fano variety''). That is to say there is an algebraic torus action on $X_{\infty}$ an element $\xi$ in the Lie algebra of this torus and a metric on the regular part of $X_{\infty}$ (locally, in $X_{\infty}$,  given by bounded potential) and which satisfies the equation (12) with respect to the holomorphic vector field generated by $\xi$.

In the case when the limiting metric is weak K\"ahler-Einstein but $X_{\infty}$ is not isomorphic to $X$ the same arguments as in 4.1 above, using the reductivity of the automorphism group, show that there is a test configuration for $X$ with central fibre $X_{\infty}$ and Futaki invariant zero. When the limit is a genuine Kahler-Ricci soliton the statment of \cite{kn:CSW} is slightly more subtle. They show that there is a destabilising test configuration for $X$, with central fibre $\overline{X}$ and strictly negative Futaki invariant,  and a further degeneration (which might be trivial) of $\overline{X}$ to $X_{\infty}$. The upshot then is
a trichotomy:
\begin{itemize} \item $X$ is K-stable: the Ricci flow converges to a Kahler-Einstein metric on $X$;
\item $X$ is K-semistable but not $K$-stable, the limit is K\"ahler-Einstein (possibly singular) but $X_{\infty}$ is not isomorphic to $X$;
\item $X$ is not K-semistable: the limit is a genuine Ricci soliton and $X_{\infty}$ is not isomorphic to $X$.
\end{itemize}

The foundation for these results is the subsequential convergence of the flow established by Chen and Wang in the earlier paper \cite{kn:CW}. That is, they prove that for any sequence $t_{i}\rightarrow \infty$ there is a subsequence $i'$ such that
$(X, \omega_{t_{i'}})$ converges (in the same sense as in the previous subsections) to a weak K\"ahler-Ricci soliton metric on a Q-Fano variety. A fundamental  difficulty in proving this is that it is not known that the Ricci curvature is bounded, either above or below, along the flow. This prevents a direct application of the Cheeger-Colding convergence theory. Results from the previous literature on K\"ahler-Ricci flow, including those of  Perelman  elaborated and extended by Sesum and Tian \cite{kn:ST}, yield three important pieces of information.
\begin{enumerate}\item The scalar curvature $R_{t}$ is bounded along the flow: $\vert R_{t}\vert\leq C_{1}$
\item The potential $h_{t}$ is bounded along the flow: $\vert h_{t}\vert\leq C_{2}$
\item There is a uniform Sobolev inequality \cite{kn:Ye}, \cite{kn:Zhang}:
$$  \Vert f \Vert_{L^{2n/n-1}}\leq C_{3}\left( \Vert \nabla f \Vert_{L^{2}} + \Vert f \Vert_{L^{2}}\right) . $$

\end{enumerate}

The general strategy of Chen and Wang is to establish a compactness theorem for segments of the flow over a fixed time interval,  say $(T-1\leq t\leq T+1)$. That is, they show that if $T_{i}$ is a sequence tending to infinity then, after passing to a subsequence, these segments of the flow converge (to a possibly singular limit). They show that the limit is \lq\lq stationary'', in that it is the solution of the K\"ahler-Ricci flow given by a K\"ahler-Ricci soliton, evolving by the action of holomorphic automorphisms. The  proofs of Chen and Wang are based on a blow-up argument and the comparison with suitable \lq\lq canonical neighbourhoods''. For $\kappa>0$, Chen and Wang define a class $\canon$ of non-compact length spaces. A space $W$ in $\canon$ is a smooth Ricci-flat K\"ahler $n$-manifold outside a closed singular set of codimension $>3$,  and with asymptotic volume ratio $\geq \kappa$,
that is $$  {\rm liminf}_{r\rightarrow \infty} v_{p}(r)\geq \kappa . $$

In their application $\kappa=\kappa(C_{3})$ is determined by the constant $C_{3}$ in the Sobolev bound. They prove from their definition that spaces in $\canon$ have many of the properties of the limit spaces treated by the Cheeger-Colding theory. In particular, they adapt that theory to show the existence of metric tangent cones. They also establish a compactness property of $\canon$, with respect to based convergence over bounded sets. This compactness means that spaces in $\canon$ satisfy certain uniform estimates. 

Chen and Wang's blow-up argument is governed by a canonical radius $cr(p)\in (0,\infty]$ which they define for any point $p$ in a Riemannian manifold $M$. This notion is in the same order of ideas as others in the literature such as the harmonic radius and curvature scale, but Chen and Wang's definition is tailored to the particular case at hand. The general idea is that $cr(p)\geq r$ if the $r$-ball centred at $p$, scaled to unit size, satisfies various definite estimates. The parameters in these estimates are chosen in line with the uniform estimates established in $\canon$, so that roughly speaking $cr(p)=\infty$ for a point in a space $W$ in $\canon$ (or more precisely $cr(p)$ is arbitarily large for a point in a Riemannian manifold which is sufficiently close to some $W\in \canon$). Now Chen and Wang establish a lower bound $cr(p)\geq \epsilon>0$ for any $p$ in a manifold $(X,\omega_{t})$ along the Ricci flow. In outline, the argument is to suppose not,  so there is a sequence of times $t_{i}$ and points $p_{i}\in (X,\omega_{t_{i}})$ such that $cr(p_{i})\rightarrow 0$. Rescaling by $r_{i}^{-1}$ they arrive at a sequence of based manifolds $(p_{i}, M_{i})$ with $cr$ bounded below by $1$ and with $cr(p_{i})=1$. They show that these manifolds converge to a space in $\canon$ and derive a contradiction from the fact that $cr=\infty$ on $\canon$.

In this argument the bound $\vert R\vert\leq C_{1}$ on the scalar curvature enters in the following way. Under the Ricci flow the scalar curvature $R$ evolves by
$$     \frac{\partial R}{\partial t}= \Delta R + \vert {\rm Ric}\vert^{2} - n $$
 After rescaling a portion of the Ricci flow, by a large factor $r_{i}^{-1}$ in the space direction and by $r_{i}^{-1/2}$ in the time direction, the scalar curvature $R'$ satisfies

    $$ \frac{\partial R'}{\partial t}= \Delta R' + \vert {\rm Ric'}\vert^{2}
- n r_{i}^{2} $$
and $\vert R'\vert \leq C_{1} r_{i}^{2}$. Thus on any region where the rescaled flows,  with a sequence of scalings $r_{i}^{-1}\rightarrow \infty$, converge in $C^{\infty}$ the limit must be a stationary Ricci-flat manifold. This is the fundamental mechanism which leads to the singular Ricci-flat spaces in $\canon$. 
The parameter $\kappa=\kappa(C_{3})$ is determined by a standard relation between the Sobolev constant and volume ratio. (If a space has small asymptotic volume ratio one can write down a compactly supported function $f$ with $\Vert\nabla f\Vert_{L^{2}}$ small compared with $\Vert f \Vert_{L^{2n/n-1}}$. If the volume ratio is less than $\kappa(C_{3})$,  such a space cannot arise as a blow-up limit of manifolds with Sobolev bound $C_{3}$.)

The $L^{2}$   construction of holomorphic sections features in two ways in Chen and Wang's arguments. One is global, to produce a projective embedding of the limit space. The tangent cone information from the blow-up limit in $\canon$ is transferred to the manifolds in the limiting sequence and the  techniques outlined in Section 3 apply. The other is local, to produce local holomorphic co-ordinates as ratios of suitable holomorphic sections. The bound on the potential $\vert h_{t}\vert \leq C_{2}$ is important here. First, since there is no lower bound on the Ricci curvature the argument based on the formula (7)does not immediately apply.  Changing the metric on the line bundle by a factor $e^{h_{t}}$ introduces an extra term which precisely cancels the Ricci curvature contribution in (7) and the bound on the $h_{t}$ means that this change does not substantially effect the estimates. Second, the evolution equation (13) gives control of the change in the metric on the line bundle in time,  and Chen and Wang are able to use this to obtain local holomorphic co-ordinates that are adapted to the metrics $\omega_{t}$ over definite time intervals. Alongside these complex geometry arguments they also use the the Ricci flow techniques of Perelman. 

The existence of a limit of the K\"ahler-Ricci flow, as established by Chen, Sun and Wang in \cite{kn:CSW} introduces further  new ideas. The Chen and Wang result  leaves open the possibility that  different sequences of times $t_{i}\rightarrow \infty$ could lead to different limits. Let ${\cal C}$ be the set of all limits that arise. By   general principles this set is connected.  One major step is to show that all $X_{\infty}$ in ${\cal C}$ can be embedded in projective space in such a way that the soliton vector fields are generated by the same fixed 1-parameter subgroup. This uses an algebro-geometric characterisation of the vector field of a Ricci soliton, via a generalisation of the Futaki invariant theory,  which leads to a rigidity property.

\subsection{ Proof by variational method}


The result proved by Berman, Boucksom and Jonsson in \cite{kn:BBJ} involves a notion of \lq\lq uniform K-stability''. Let ${\cal X}$ be a test configuration for a polarised manifold $(X,L)$, so we have a $\bC^{*}$-action on the central fibre $X_{0}$. Suppose first that $X_{0}$ is smooth and we fix an $S^{1}$-invariant K\"ahler metric in the class $c_{1}({\cal L}_{0})$ which yields a symplectic structure. The $S^{1}$-action is generated by a Hamiltonian function $H$ which we can normalise to have maximum value $0$. Then we define    $\Vert {\cal X}\Vert$ to be the $L^{1}$ norm of $H$. This is a quantity which is independent of the choice of metric in the cohomology class.  The definition can be extended to any scheme, using the asymptotics of the trace of the action on sections of ${\cal L}^{k}$ as $k\rightarrow\infty$, similar to  the definition of the Futaki invariant. Then $(X,L)$ is said to be {\it uniformly K-stable} if there is some $\epsilon>0$ such that
\begin{equation}     {\rm Fut}({\cal X})\geq \epsilon \Vert {\cal X}\Vert \end{equation}for all non-trivial test configurations ${\cal X}$.  The main result of \cite{kn:BBJ} is that for a Fano manifold $X$, polarised by $K_{X}^{-1}$ and  with finite automorphism group,  the existence of a K\"ahler-Einstein metric is equivalent to uniform $K$-stability. Note that {\it a priori} uniform $K$-stability  is a stronger condition than $K$-stability although {\it a posteriori} they are equivalent, for Fano manifolds with finite automorphism group. One can consider many other  norms on test configurations and the general notion of uniform stability goes back to the thesis of Szekelyhidi. It is also related to another  variant of  $K$-stability developed by Szekelyhidi which considers filtrations of the co-ordinate ring $\bigoplus H^{0}(X,L^{k})$ \cite{kn:Gabor1}.  A definition of uniform stability which turns out to be equivalent to that in \cite{kn:BBJ} was given by Dervan \cite{kn:Dervan}.

To give some indication of the proof in \cite{kn:BBJ}  we begin by considering an analogous situation in finite dimensions. Let $V$ be a complete (finite-dimensional) Riemannian manifold with the property that the each two points can be joined by a unique geodesic segment---for example a simply connected manifold of non-positive curvature. Let $v_{0}\in V$ be a base point and let $F$ be  function on $V$ which is convex along geodesics. If $\gamma:[0,\infty)\rightarrow V$ is a geodesic ray emanating from $v_{0}$, parametrised by arc length,  the ratio $\gamma(t)/t$ is increasing with $t$ and we can define the {\it asymptotic slope} $ S_{\gamma}\in [0,\infty]$ to be the limit as $t\rightarrow \infty$. If $S_{\gamma}>0$ then for any $\delta< S_{\gamma}$ there is a $C_{\gamma}$ such that \begin{equation}   F(\gamma(s)) \geq \delta s -C_{\gamma}. \end{equation}
 An elementary argument, hinging on the compactness of the set of geodesics through $v_{0}$,  shows that if there is some $\delta>0$ such that if $S_{\gamma}\geq \delta$ for all such rays $\gamma$ then given some  $\delta'<\delta$ we can find $C$ such  
\begin{equation}  F(v)\geq \delta' d(v,v_{0}) - C  \end{equation}
for all $v\in V$. It also follows easily that $F$ attains a minimum in $V$.

The relevance of this to the K\"ahler-Einstein problem on a Fano manifold $X$  is that,  as we have discussed in Section 2, a K\"ahler-Einstein metric can be seen as a critical point of the Mabuchi functional ${\cal F}$ on the space of K\"ahler metrics ${\cal H}$.  The pair $({\cal H}, {\cal F})$ has many properties analogous to $(V,F)$ above. There is a Mabuchi metric which makes ${\cal H}$ formally a symmetric space of non-positive curvature and ${\cal F}$ is convex along geodesics. The programme roughly speaking, is to extend the arguments above to this infinite-dimensional setting and to relate the asymptotics of the Mabuchi functional, analogous to the asymptotic slope $S_{\gamma}$,  to the condition (17) on test configurations.

A point $\omega$ in ${\cal H}$ defines a volume form on $X$ and the tangent space to ${\cal H}$ at $\omega$ can be identified with the functions $\delta\phi$ of integral $0$. For each $p\geq1$ the $L^{p}$ norm defines a Finsler structure on ${\cal H}$. The case $p=2$ gives the infinite dimensional Riemannian structure first considered by Mabuchi \cite{kn:Mab1} but the case $p=1$ is also important, as shown by Darvas \cite{kn:Darvas}.
The completion, in the metric defined by this Finsler structure, is a space ${\cal E}^{1}$ of currents defined by \lq\lq finite energy'' potentials.  Geodesics in ${\cal H}$ also have a good geometric meaning. Smooth geodesics segments parametrised by an interval $[a,b]$ correspond to $S^{1}$-invariant closed $(1,1)$-forms $\Omega$ on the product $X\times S^{1}\times [a, b]$ which satisfy $\Omega^{n+1}=0$ and which restrict to a metric in ${\cal H}$ on each copy of $X$ in the product. (The different Finsler structures share the same geodesics.) If we stay in the smooth category it is not true that any two points in ${\cal H}$ can be joined by a geodesic \cite{kn:Lempert} but this is the case if one relaxes the definitions to allow singularities---in fact forms with $C^{1,1}$ potentials---as shown by Chen \cite{kn:XC}. It is an elementary calculation that the Mabuchi functional is convex along smooth geodesics. The convexity in the general case is a deep recent result of Berman and Berndtsson \cite{kn:BB}.

Now  fix a base point $\omega_{0}$ in ${\cal H}$---analogous to $v_{0}\in V$.  According to Darvas, the distance to $\omega_{0}$ in the $L^{1}$-Finsler structure is equivalent to functional $J$ which is well-known in the literature and  which can be characterised by the property that $J(\omega_{0})=0$ and
$$  \delta J = \int_{X} \delta \phi (\omega^{n}-\omega_{0}^{n}) . $$
Thus the analogue of (19) is an inequality.
\begin{equation}    {\cal F}(\omega)\geq \delta' J(\omega)-C\end{equation}
for some fixed $\delta'>0, C$ and $\omega$. This is sometimes referred to as the \lq\lq properness'' of the Mabuchi functional. Tian showed \cite{kn:Ti}  that if this inequality (20) holds then a K\"ahler-Einstein metric exists. The Mabuchi functional is decreasing along the continuity path,  as discussed in 4.2 above, and so ${\cal F}(\omega_{s})$ controls $J(\omega_{s})$ and Tian showed that  this allows the continuity path to be continued to $s=1$. Darvas and Rubinstein have given another  proof of this result, and generalisations, recently \cite{kn:DR}.

There is a large circle of results relating geodesic rays in ${\cal H}$   to test configurations, see \cite{kn:PS} for example.  Using the conformal equivalence between $S^{1}\times [0,\infty)$ and the punctured disc $\Delta^{*}$, geodesic rays correspond to $S^{1}$-invariant closed $(1,1)$-forms $\Omega$ on the product $X\times \Delta^{*}$ with $\Omega^{n+1}=0$ and which are positive on each $X\times \{t\}$. For certain purposes  one can work with {\it subgeodesics} which correspond to positive definite forms $\Omega$. In particular if ${\cal X}$ is a test configuration we can consider a \lq\lq smooth'' metric  $\Omega$ on ${\cal X}$. For example we can embed ${\cal X}$ in some  $\Delta\times \bC\bP^{N}$ and take the restriction of a metric on the ambient manifold. The ${\bf C}^{*}$-action on ${\cal X}$ gives an open embedding $X\times \Delta^{*}$ and we get a subgeodesic ray $\omega_{s}$  in ${\cal H}$.   Boucksom, Hisamoto and Jonsson \cite{kn:BHJ} prove that (if the central fibre is reduced)
 $$     J(\omega_{s}) \sim \Vert {\cal X}\Vert s\ \ \  ;\ \ \ \  {\cal F}(\omega_{s})\sim {\rm Fut}({\cal X}) s ,  $$
as $s\rightarrow \infty$. 
(There are many earlier results of this kind in the literature, under various hypotheses.) 
Thus the uniform stability condition is equivalent to the statement that there is a $\delta>0$ such that for any such subgeodesic ray,  arising from a test configuration,  we have
\begin{equation}  {\cal F}(\omega_{s})\geq \delta J(\omega_{s})- C, \end{equation}
where $C$ depends on the ray. 

There are now two main aspects to the proof. 
(The exposition in \cite{kn:BBJ} involves some sophisticated
techniques which go well beyond this writers knowledge, and indeed \cite{kn:BBJ}
is described by the authors as an outline to be followed by a more detailed
version. What we write below is extremely sketchy.)
\begin{itemize} \item 
 To pass from the  subgeodesic rays arising from test configurations to general geodesic rays and establish an inequality (21) along any geoedesic ray.

A sub-geodesic ray comes with a family of K\"ahler potentials which can be viewed as a metric on the pull-back of $L$ to $X\times \Delta^{*}$ or as a singular metric on the pull-back to $X\times \Delta$. This metric defines  a multiplier ideal sheaf: the local holomorphic sections which are in $L^{2}$ with respect to the metric. Fundamental results of Nadel show that this is  a  coherent sheaf and so one can construct  corresponding blow-ups of $X\times \Delta$ along the powers of this ideal sheaf,  which yield  test configurations. Bermann, Boucksom and Jonnson  use these to approximate the original ray by those arising from test configurations and eventually to pass from the algebro-geometric uniform K-stability hypothesis to the estimate  (18) on general geodesic rays. For these purposes they also work with the Ding functional (which we mentioned briefly in Section 2).  They also use ideas and results from  non-Archimedean geometry.

\item   To carry the elementary arguments from the finite dimensional model
over to the infinite dimensional situation.
Results of Berman, Boucksom, Eyssidieux, Guedj and Zeriahi \cite{kn:BBEGZ} are used here to give the relevant compactness
property in ${\cal E}^{1}$ for geodesics segments with bounded Mabuchi functional. \end{itemize}

   These variational techniques based on convex geometry in the space ${\cal H}$ of K\"ahler metrics have been used by Darvas and Rubinstein \cite{kn:DR}  and Berman, Darvas and Lu \cite{kn:BDL} to produce interesting results in the more general framework of constant scalar curvature metrics. The outstanding problem is to prove the regularity of weak solutions produced by minimising  the Mabuchi functional.



\end{document}